\documentclass[12pt]{article} 
\usepackage{amsmath, amsthm}
\begin{document}
\input amssym.def
\input amssym

\setcounter{equation}{0}

\newcommand{\wt}{\mbox{wt}}
\newcommand{\spa}{\mbox{span}}
\newcommand{\Res}{\mbox{Res}}
\newcommand{\End}{\mbox{End}}
\newcommand{\Ind}{\mbox{Ind}}
\newcommand{\Hom}{\mbox{Hom}}
\newcommand{\Mod}{\mbox{Mod}}
\newcommand{\m}{\mbox{mod}\ }
\renewcommand{\theequation}{\thesection.\arabic{equation}}

\def \Aut{{\rm Aut}}
\def \Z{\Bbb Z}
\def \M{\Bbb M}
\def \C{\Bbb C}
\def \R{\Bbb R}
\def \Q{\Bbb Q}
\def \N{\Bbb N}
\def \ann{{\rm Ann}}
\def \<{\langle} 
\def \o{\omega}
\def \O{\Omega}
\def \M{{\cal M}}
\def \1t{\frac{1}{T}}
\def \>{\rangle} 
\def \t{\tau }
\def \a{\alpha }
\def \e{\epsilon }
\def \l{\lambda }
\def \L{\Lambda }
\def \g{\gamma}
\def \b{\beta }
\def \om{\omega }
\def \o{\omega }
\def \cg{\chi_g}
\def \ag{\alpha_g}
\def \ah{\alpha_h}
\def \ph{\psi_h}
\def \nor{\vartriangleleft}
\def \V{V^{\natural}}
\def \voa{vertex operator algebra\ }
\def \v{vertex operator algebra\ }
\def \1{{\bf 1}}
\def \be{\begin{equation}\label}
\def \ee{\end{equation}}
\def \pf {\noindent {\bf Proof:} \,}
\def \bl{\begin{lem}\label}
\def \el{\end{lem}}
\def \ba{\begin{array}}
\def \ea{\end{array}}
\def \bt{\begin{thm}\label}
\def \et{\end{thm}}
\def \ch{{\rm ch}}
\def \br{\begin{rem}\label}
\def \er{\end{rem}}
\def \ed{\end{de}}
\def \bp{\begin{prop}\label}
\def \ep{\end{prop}}

\newtheorem{th1}{Theorem}
\newtheorem{ree}[th1]{Remark}
\newtheorem{thm}{Theorem}[section]
\newtheorem{prop}[thm]{Proposition}
\newtheorem{coro}[thm]{Corollary}
\newtheorem{lem}[thm]{Lemma}
\newtheorem{rem}[thm]{Remark}
\newtheorem{de}[thm]{Definition}
\newtheorem{hy}[thm]{Hypothesis}

\begin{center}
{\Large {\bf Twisted representations of  vertex operator superalgebras}}\\
\vspace{0.5cm}
 
Chongying Dong\footnote{Supported by NSF grants, a China NSF grant
and a faculty research fund from the University of 
California at 
Santa Cruz.} and 
Zhongping Zhao
\\
Department of Mathematics, University of
California, Santa Cruz, CA 95064
\hspace{1.5 cm}
\end{center}

\begin{abstract} This paper gives an analogue of $A_g(V)$ theory
for a vertex operator superalgebra $V$ and an automorphism 
$g$ of finite order. The relation between the $g$-twisted $V$-modules and 
$A_g(V)$-modules is established. It is proved that if $V$ is $g$-rational,
then $A_g(V)$ is finite dimensional semisimple associative algebra and
there are only finitely many irreducible $g$-twisted $V$-modules. 
\end{abstract}

\section{Introduction}
\setcounter{equation}{0}

The twisted sectors or twisted modules are basic ingredients in orbifold 
conformal field theory (cf. \cite{FLM1}, \cite{FLM2}, \cite{FLM3}, \cite{Le1},
\cite{Le2}, \cite{DHVW}, \cite{DVVV}, \cite{DL2}, \cite{DLM2}).
 The notion of twisted module \cite{FFR},\cite{D}
is derived from the properties of twisted vertex operators  
for finite automorphisms of even lattice vertex operator algebras 
constructed in \cite{Le1}, \cite{Le2} and \cite{FLM2}, also see \cite{DL2}. 
In this paper we study the twisted modules for an arbitrary vertex
operator superalgebra following \cite{Z}, \cite{KW} and \cite{DLM2}.

An associative algebra $A(V)$ was introduced in \cite{Z} for every vertex
operator algebra $V$ to study the representation theory for vertex operator
algebra. The main idea is to reduce the study of representation
theory for a vertex operator algebra to the study of representation theory for
an associative algebra. This approach has been very successful and the 
irreducible modules for many well-known vertex operator algebras
have been classified by using the associative algebras. This theory has
been extended to the vertex operator superalgebras in \cite{KW}
and has been further generalized to the twisted representations
for a vertex operator algebra in \cite{DLM2}.

This paper is a ``super analogue'' of \cite{DLM2}. We construct 
an associative algebra $A_g(V)$ for any vertex operator superalgebra
$V$ together with an automorphism $g$ of finite order. Then
the vacuum space of any admissible $g$-twisted $V$-module 
becomes a module for $A_g(V).$  On the other hand one can construct
a `universal' admissible $g$-twisted $V$-module from any
$A_g(V)$-module. This leads to a one to one correspondence 
between the set of inequivalent admissible $g$-twisted $V$-modules
and the set of simple $A_g(V)$-modules. As in the case of vertex
operator algebra, if $V$ is $g$-rational then $A_g(V)$ is 
a finite dimensional semisimple associative algebra. 

The ideas of this paper and other related papers are very
natural and go back to the theory 
of highest weight modules for Kac-Moody Lie algebras and other
Lie algebras with triangular decompositions. In the classical 
highest weight module theory, the highest weight or highest weight vector
determines the highest weight module structure to some extend (different
highest weight modules can have the same highest weight). The role of the
vacuum space for an admissible twisted module is similar to the
role of the highest weight space in a highest weight module. So from this
point of view, the $A_g(V)$ theory is a natural extension
of highest weight module theory in the representation theory of vertex
operator superalgebras.

A  vertex operator superalgebra has a canonical automorphism $\sigma$ 
of order $2$ arising from the structure of superspace.
The $\sigma$-twisted modules which are called the Ramond sector in the
literature play very important roles in the study of geometry.
Important topological invariants such as elliptic genus and certain
Witten genus can be understood as graded trace functions on the Ramond
sectors constructed from the manifolds. 
It is expected that the theory developed in this paper will have  
applications in geometry and physics.
 
Since the setting and most results in this paper are similar to those in
\cite{DLM2} we only provide the arguments which are either new or need
a lot of modifications. We refer the reader to \cite{DLM2} for details. 

The organization of this paper is similar to that of \cite{DLM2}.
We review  the definition of  vertex operator superalgebra
and define various notions of $g$-twisted $V$-modules in section 2.
In section 3, we introduce 
the algebra $A_{g}(V)$ for VOSA $V$. Section 4 
is devoted to the study of Lie superalgebra $V[g]$ which is kind of
twisted affinization of $V.$  A $weak \ \ g$-twisted $V$-module is naturally a
$V[g]$-module. In section 5,
we construct the functor $\O$ which sends a weak $g$-twisted $V$-module to
an $A_g(V)$-module. We construct another functor  $L$
from the category of $A_g(V)$-modules to the category of
admissible $g$-twisted $V$-modules in Section 6. That is, for any
$A_g(V)$-module $U$ we can construct a kind of
``generalized Verma module'' $\bar M(U)$ which is the universal admissible 
$g$-twisted $V$-module generated by $U.$ It is proved that
there is a 1-1 correspondence between the irreducible objects
in these two categories. Moreover if $V$ is $g$-rational, then
$A_g(V)$ is a finite dimensional semisimple associative
algebra. We discuss some examples of vertex operator superalgebras
constructed from the free fermions and their twisted modules in
Section 7.

\section{Vertex Operator superalgebra and twisted modules}
\setcounter{equation}{0}

We review the definition of vertex operator superalgebra (cf. \cite{B},
\cite{FLM3}, \cite{DL1}) and
various notions of twisted modules in this section (cf. \cite{D}, \cite{DLM2},
\cite{FFR}, \cite{FLM3}, \cite{Z}). 

Recall that a super vector space  is a $\Bbb Z_{2}$-graded vector space
 $V=V_{\bar{0}}\oplus V_{\bar{1}}$. The elements in $V_{\bar{0}}$ 
(resp. $V_{\bar{1}}$) are called even (resp. odd). Let $\tilde{v}$ 
be $0$ if $v\in V_{\bar{0}}$, and $1$ if  $v\in V_{\bar{1}}$.

\begin{de} {\rm A  vertex operator superalgebra is a
$\frac{1}{2}\Bbb Z_{+}$-graded super vector space
\begin{equation}\label{g2.1}
V=\bigoplus_{n\in{ \frac{1}{2}\Bbb Z_{+}}}V_n= V_{\bar{0}}\oplus V_{\bar{1}}.
\end{equation}
with  $V_{\bar{0}}=\sum_{n\in\Z}V_n$ and 
$V_{\bar{1}}=\sum_{n\in\frac{1}{2}+\Z}V_n$ 
satisfying $\dim V_{n}< \infty$ for all $n$ and $V_m=0$ if $m$ is sufficiently
small.  $V$ is   equipped with a linear map
 \begin{align}\label{0a3}
& V \to (\mbox{End}\,V)[[z,z^{-1}]] ,\\
& v\mapsto Y(v,z)=\sum_{n\in{\Z}}v_nz^{-n-1}\ \ \ \  (v_n\in
\mbox{End}\,V).\nonumber
\end{align}
and with two distinguished vectors ${\bf 1}\in V_0,$ $\omega\in
V_2$ satisfying the following conditions for $u, v \in V,$ and $m,n\in\Z:$
\begin{align} \label{0a4}
& u_nv=0\ \ \ \ \ {\rm for}\ \  n\ \ {\rm sufficiently\ large};  \\
& Y({\bf 1},z)=Id_{V};  \\
& Y(v,z){\bf 1}\in V[[z]]\ \ \ {\rm and}\ \ \ \lim_{z\to
0}Y(v,z){\bf 1}=v;\\
& [L(m),L(n)]=(m-n)L(m+n)+\frac{1}{12}(m^3-m)\delta_{m+n,0}c ;\\
& \frac{d}{dz}Y(v,z)=Y(L(-1)v,z);\\
& L(0)|_{V_n}=n.
\end{align}
where $L(m)=\o_{ m+1}, $ that is,
$$Y(\o,z)=\sum_{n\in\Z}L(n)z^{-n-2};$$
and the {\em Jacobi identity} holds: 
\begin{equation}\label{2.8}
\begin{array}{c}
\displaystyle{z^{-1}_0\delta\left(\frac{z_1-z_2}{z_0}\right)
Y(u,z_1)Y(v,z_2)-(-1)^{\tilde{u}\tilde {v}}z^{-1}_0\delta\left(\frac{z_2-z_1}{-z_0}\right)
Y(v,z_2)Y(u,z_1)}\\
\displaystyle{=z_2^{-1}\delta
\left(\frac{z_1-z_0}{z_2}\right)
Y(Y(u,z_0)v,z_2)}.
\end{array}
\end{equation}
where 
$\delta(z)=\sum_{n\in {\Bbb Z}}z^n$ and  $(z_i-z_j)^n$ is 
expanded as a formal power series in $z_j.$ Throughout the paper,
$z_0,z_1,z_2,$ etc. are independent commuting formal variables.}
\end{de}

Such a vertex operator superalgebra may be denoted by  $V=(V,Y,{\bf 1},\omega).$
In the case $V_{\bar 1}=0,$ this is exactly the definition
of vertex operator algebra given in [FLM3].

\begin{de}{\rm Let $V$ be a vertex operator superalgebra. 
An {\em automorphism} $g$ of $V$ is a
linear automorphism of $V$ preserving  $\omega$ such that
the actions of $g$ and $Y(v,z)$ on $V$ are compatible in the sense
that $$gY(v,z)g^{-1}=Y(gv,z)$$ for $v\in V.$}
\end{de}

Note that any automorphism of $V$ commutes with $L(0)$ and preserves
each homogeneous space $V_n.$ As a result, any automorphism preserves
$V_{\bar 0}$ and $V_{\bar 1}.$ 

Let $\Aut(V)$ be the group of automorphisms of $V.$ There is a
special automorphism  $\sigma \in \Aut(V)$ such that $\sigma|V_{\bar 0}=1$
and $ \sigma|V_{\bar 1}=-1.$ It is clear that $\sigma$ is a central element
of $\Aut(V).$ 

Fix $g\in \Aut(V)$ of order $T_0.$ Let $o(g\sigma)=T.$  Denote the decompositions of $V$ into eigenspaces with respect to
the actions of $g\sigma$  and $g$ as follows
\begin{eqnarray}
& &V=\oplus_{r\in \Z/T\Z}V^{r*} \label{g2.9}\\
& & V=\oplus_{r\in \Z/T_{0}\Z}V^{r} \label{g2.10} 				\end{eqnarray}
where $V^{r*}=\{v\in V|g\sigma v=e^{2\pi ir/T}v\}$ and $V^{r}=\{v\in V|gv=e^{2\pi ir/T_0}v\}$

\begin{de} {\rm A weak $g$-twisted $V$-module $M$ is a vector space equipped 
with a linear map
$$\begin{array}{l}
V\to (\End\,M)[[z^{1/T_0}, z^{-1/T_0}]\\
v\mapsto\displaystyle{ Y_M(v,z)=\sum_{n\in\frac{1}{T_0}\Z}v_nz^{-n-1}\ \ \ (v_n\in
\End\,M)}
\end{array}$$
which satisfies that for all $0\leq r\leq T_0-1,$ $u\in V^r,$ $v\in V,$ 
$w\in M,$
\begin{eqnarray}\label{g2.11}
& &Y_M(u,z)=\sum_{n\in \frac{r}{T_0}+\Z}u_nz^{-n-1} \label{1/2} ;\\ 
& &u_lw=0 \ \ \  				
\mbox{for}\ \ \ l>>0\label{vlw0};\\
& &Y_M(\1,z)=Id_{M};\label{vacuum}
\end{eqnarray}
 \begin{equation}\label{2.14}
\begin{array}{c}
\displaystyle{z^{-1}_0\delta\left(\frac{z_1-z_2}{z_0}\right)
Y_M(u,z_1)Y_M(v,z_2)-(-1)^{\tilde{u}\tilde{v}}z^{-1}_0\delta\left(\frac{z_2-z_1}{-z_0}\right)
Y_M(v,z_2)Y_M(u,z_1)}\\
\displaystyle{=z_2^{-1}\left(\frac{z_1-z_0}{z_2}\right)^{-r/T_0}
\delta\left(\frac{z_1-z_0}{z_2}\right)
Y_M(Y(u,z_0)v,z_2)}.
\end{array}
\end{equation}}
\end{de}

Following the arguments in  [DL1] one can prove
that the twisted Jacobi identity is equivalent to the following
associativity formula 
\begin{eqnarray}\label{ea2.15}
(z_{0}+z_{2})^{k+\frac{r}{T_0}}Y_{M}(u,z_{0}+z_{2})Y_{M}(v,z_{2})w
=(z_{2}+z_{0})^{k+\frac{r}{T_0}}Y_M(Y(u,z_0)v,z_2)w.
\end{eqnarray}
where $w\in M$ and $k\in\Bbb Z_{+}$ s.t $z^{k+\frac{r}{T_0}}Y_{M}(u,z)w$ involves only nonnegative integral powers of $z,$ and 
commutator relation
\begin{eqnarray}\label{g2.16}
& &\ \ \ \  [Y_{M}(u,z_{1}),Y_{M}(v,z_{2})]\nonumber\\
& &=\Res_{z_{0}}z_2^{-1}\left(\frac{z_1-z_0}{z_2}\right)^{-r/T_0}
\delta\left(\frac{z_1-z_0}{z_2}\right)Y_M(Y(u,z_0)v,z_2).\label{ec}
\end{eqnarray}

We need the following Lemma later.
\begin{lem}\label{ladd}
The associativity formula (\ref{ea2.15}) is equivalent to the following:
\begin{eqnarray*}
(z_{0}+z_{2})^{m+\frac{s}{T}}Y_{M}(u,z_{0}+z_{2})Y_{M}(v,z_{2})w
=(z_{2}+z_{0})^{m+\frac{s}{T}}Y_M(Y(u,z_0)v,z_2)w
\end{eqnarray*}
for $u\in V^{s*}$ and some $m\in\frac{1}{2}\Z_+$ with $m=\wt\,u$ modulo
$\Z$ such that  $z^{m+\frac{s}{T}}Y_{M}(u,z)w$ involves only nonnegative integral powers of $z.$
\end{lem}

\pf Let $u\in V^r.$ It is enough to prove that $\wt\,u+\frac{s}{T}$ and $
\frac{r}{T_0}$ are congruent modulo $\Z.$ It is easy to see that
$s\equiv \frac{T}{2}\tilde{u}+r$ modulo $\Z$  if $T_0$ is even and 
$s\equiv \frac{T}{2}\tilde{u}+2r$ modulo $\Z$ if $T_0$ is odd.  Thus
$\wt\,u+\frac{s}{T}$ is congruent to $\wt\,u+\frac{1}{2}\tilde{u}+\frac{r}{T_0}.$
Since $\frac{1}{2}\tilde{u}$ and $\wt\,u$ are congruent modulo $\Z,$ the result
follows immediately. \qed

Equating the coefficients of
$z_1^{-m-1}z_2^{-n-1}$ in (\ref{ec}) yields
\begin{eqnarray}\label{g2.17}
[u_{m},v_{n}]=\sum_{i=0}^{\infty}
\left(\begin{array}{c}m\\i\end{array}\right)(u_{i}v)_{m+n-i}.
\end{eqnarray}

We may also deduce from
(\ref{1/2})-(\ref{2.14}) the usual Virasoro algebra axioms, namely that
if $Y_M(\o,z)=\sum_{n\in\Z}L(n)z^{-n-2}$ then
\begin{equation}\label{g2.18}
[L(m),L(n)]=(m-n)L(m+n)+\frac{1}{12}(m^3-m)\delta_{m+n,0}c
\end{equation}
and 
\begin{equation}\label{g2.19}
\frac{d}{dz}Y_M(v,z)=Y_M(L(-1)v,z)
\end{equation}
(cf. [DLM1]). 

The homomorphism and isomorphism of weak twisted modules are defined in an
obvious way. 

\begin{de} An {\em admissible} $g$-twisted $V$-module 
is a  weak $g$-twisted $V$-module $M$ which carries a 
$\frac{1}{T}{\Z}_{+}$-grading 
\begin{equation}\label{g2.22}
M=\oplus_{n\in\frac{1}{T}\Z_+}M(n)
\end{equation}
satisfying 
\begin{eqnarray}\label{g2.23}
v_{m}M(n)\subseteq M(n+\wt v-m-1)
\end{eqnarray}
for homogeneous $v\in V.$
\ed

\begin{de}\label{d2.4} An ordinary $g$-{\em twisted $V$-module} is
a weak $g$-twisted $V$-module
\begin{equation}\label{g2.21}
M=\coprod_{\lambda \in{\C}}M_{\lambda} 
\end{equation}
such that $\dim M_{\l}$ is finite and for fixed $\l,$ $M_{\frac{n}{T}+\l}=0$
for all small enough integers $n$ where
 $M_{\l}=\{w\in M|L(0)w=\l w\}.$
\end{de}

The admissible $g$-twisted $V$-modules form a subcategory 
of the weak $g$-twisted $V$-modules. It is easy to prove that
an ordinary $g$-twisted $V$-module is admissible.
Shifting the grading of
an admissible $g$-twisted module gives an isomorphic admissible $g$-twisted
$V$-module. A {\em simple}
object in this category is an admissible $g$-twisted $V$-module $M$ such that
$0$ and $M$ are the only graded submodules.

We say that $V$ is $g$-{\em rational} if every admissible $g$-twisted
$V$-module is completely reducible, i.e., a direct sum of simple admissible
$g$-twisted modules. $V$ is called {\em rational} if $V$ is $1$-rational.
$V$ is called {\em holomorphic} if $V$ is rational and $V$ is the
only irreducible $V$-module up to isomorphism.

If $M\!=\!\oplus_{n\in \frac{1}{T}{\Z}_{+}}\!M(n)$
is an admissible $g$-twisted $V$-module, the contragredient module $M'$ 
is defined as follows: 
\be{g2.25}
M'=\oplus_{n\in \frac{1}{T}{\Z}_{+}}M(n)^{*}
\end{equation}
where $M(n)^*=\Hom_{\C}(M(n),\C).$ The vertex operator 
$Y_{M'}(a,z)$ is defined for $a\in V$ via  
\begin{eqnarray}\label{g2.26}
\langle Y_{M'}(a,z)f,u\rangle= \langle f,Y_M(e^{zL(1)}(-z^{-2})^{L(0)}a,z^{-1})u\rangle
\end{eqnarray}
where $\<\cdot\>$ denotes the natural paring between $M'$ and $M.$ 
Then we have the following \cite{FHL}:
\bl{l2.9} $(M',Y_{M'})$ is an admissible $g^{-1}$-twisted $V$-module.
\el

Lemma \ref{l2.9} is needed in the proof of several results in Section
6 although we do not intend to give these proofs (cf. \cite{DLM2}).  

\section{The associative algebra $A_{g}(V)$}
\setcounter{equation}{0}

Let  $r$ be an integer between $0$ and $T-1$ ( or $T_0-1$).
We will also use $r$ to denote  its residue class modulo $T$ or
$T_0.$ 
For homogeneous $u\in V^{r*},$ we set  $\delta_{r}=1$ if $r=0$ and  $\delta_{r}=0$ if $r\ne 0.$ 
Let  $v\in V$ we  define
\begin{eqnarray}\label{g3.1}
u\circ_g v=\Res_{z}\frac{(1+z)^{{\wt u}-1+\delta_{r}+\frac{r}{T}}}{z^{1+\delta_{r}}}Y(u,z)v
\end{eqnarray}
where $(1+z)^{\alpha}$ for $\alpha\in\C$ is to
be expanded in nonnegative integer powers of $z.$ 
Let $O_g(V)$ be the linear span of all $u\circ_g v$ and define the linear space
$A_g(V)$ to be the quotient $V/O_{g}(V).$ We will use 
$A(V), O(V), u\circ v,$ when $g=1.$  The $A(V)$ was constructed in \cite{KW}
and if $V$ is a vertex operator, $A_g(V)$ was constructed in \cite{DLM2}.

\begin{lem}\label{l3.1} If $r\ne 0$ then $V^{r*}\subseteq O_{g}(V).$
\end{lem}

\pf  The proof is the same as that of Lemma 2.1 in \cite{DLM2}.
\qed

Let $I=O_g(V)\cap V^{0*}.$ Then  $A_g(V) \simeq V^{0*}/I$ (as linear spaces).
Since $O(V^{0*})\subset I,$ $A_g(V)$ is a quotient of $A( V^{0*}).$

We now define a product $*_g$ on $V$ which will induce an associative
product in $A_g(V).$ Let $r,u$ and $v$ be as
above and  set
\begin{equation}\label{3.2}
u*_gv=\left\{
\begin{array}{ll}
\Res_z(Y(u,z)\frac{(1+z)^{{\wt}\,u}}{z}v)
 & {\rm if}\ r=0\\
0  & {\rm if}\ r>0.
\end{array}\right.
\end{equation} 

As in \cite{DLM2} we extend $*_g$ linearly to obtain a bilinear product  on 
$V.$ Then the restriction of $*_g$ to $V^{0*}$ coincides with that of 
\cite{KW}. As before we will use $*$ instead of $*_g$ if $g=1.$  
If $u\in V^{0*}$ then  we can write (\ref{3.2}) as the following form
\begin{equation}\label{3.3}
u*_gv=\sum_{i=0}^{\infty}\binom{\wt u}{i}u_{i-1}v.
\end{equation}

\begin{lem}\label{l3.2} (i) Assume that $u\in V^{r*}$ is homogeneous,
$v\in V$ and $m\ge n\ge 0.$ Then 
$$\Res_z\frac{(1+z)^{{\wt}u-1+\delta_{r}+\frac{r}{T}+n}}{z^{m+\delta_{r}+1}}Y(u,z)v\in O_{g}(V).$$

(ii) Assume that $u,v\in V^{0*}$ are homogeneous. Then
$$u*v-(-1)^{\tilde {u} \tilde {v}}\Res_z\frac{(1+z)^{\wt v-1}}{z}Y(v,z)u\in O(V^{0*})$$
and

(iii) $u*v-(-1)^{\tilde {u} \tilde {v}}v*u-\Res_z(1+z)^{\wt u-1}Y(u,z)v\in O(V^{0*}).$
 \end{lem}
\pf See the proofs of Lemmas 2.1.2 and 2.1.3 of \cite{Z} by noting that
 $$Y(u,z)v\equiv(-1)^{\tilde {u} \tilde {v}}(1+z)^{-wtu-wtv}Y(v,\frac{-z}{1+z}u)\ \ mod\ \ O_{g}(V^{0*})$$
and 
$$ Y(u,z)v=(-1)^{\tilde {u} \tilde {v}}e^{zL(-1)}Y(v,-z)u$$
(cf. \cite{KW}). \qed
 
Here is our first main result.
\begin{thm}\label{t3.3} (i)  $A_g(V)$ is an associative algebra 
under $*_g.$

$(ii)$ $1+O_g(V)$ is the unit of $A_g(V).$

$(iii)$ $\o+O_g(V) $ lies in the center of $A_g(V).$
\end{thm} 
 \pf First we prove that $*_g$ is well defined on $A_g(V).$ It is equivalent
to prove that $O_g(V)$ is a 2-sided ideal of $V$ under $*_g.$  
 Since $V^{r*}*_g V=0$ if $r\neq 0$ and $V^{0*}*_g V^{r*}\subseteq V^{r*} \subseteq O_g(V)$, it is enough to prove that 
$I=O_g(V)\cap V^{0*} $ is a two sided ideal of $V^{0*}$ under $*.$  
The proof given here is similar to that of Proposition 2.3 of \cite{DLM2}.

Let $c\in V^{0*}$ be homogeneous and $u\in I$. We need to show that both
  \begin{eqnarray}\label{g3.5}
c*u=\Res_{z}\left(\frac{(1+z)^{{\wt}c}}{z}Y(c,z)u\right)
\end{eqnarray}
and
\begin{eqnarray}
u*c\equiv\Res_{z}\frac{(1+z)^{{\wt}c-1}}{z}Y(c,z)u\ (\m I)
\end{eqnarray}
lie in $I.$ 

From (\ref{g2.9})
it suffices to take $u=a\circ_g b$ where $a\in V^{r*}$ and
$b\in V^{(T-r)*}$ are both homogeneous. 
Set $x_{0}=c*u,$ $x_1=u*c$ and recall the twisted Jacobi identity on $V:$
\begin{equation}\label{g3.6}
\begin{array}{c}
\displaystyle{z^{-1}_0\delta\left(\frac{z_1-z_2}{z_0}\right)
Y(c,z_1)Y(a,z_2)b-(-1)^{\tilde {c} \tilde {a}}z^{-1}_0\delta\left(\frac{z_2-z_1}{-z_0}\right)
Y(a,z_2)Y(c,z_1)}b\\
\displaystyle{=z_2^{-1}\delta\left(\frac{z_1-z_0}{z_2}\right)
Y(Y(c,z_0)a,z_2)b}.
\end{array}
\end{equation}

For $\varepsilon=0$ or $1,$  (\ref{g3.6}) implies:
\begin{eqnarray*}
& &x_{\varepsilon}=\Res_{z_{1}}\frac{(1+z_{1})^{{\wt}c-\varepsilon}}{z_{1}}
Y(c,z_{1})\Res_{z_2}\frac{(1+z_2)^{{\wt}a-1+\delta_{r}+\frac{r}{T}}}{z_2^{1+\delta_{r}}}Y(a,z_2)b\\
& &=\Res_{z_{1}}\Res_{z_2}\frac{(1+z_{1})^{{\wt}c-\varepsilon}}{z_{1}}
Y(c,z_{1})\frac{(1+z_2)^{{\wt}a-1+\delta_{r}+\frac{r}{T}}}{z_2^{1+\delta_{r}}}Y(a,z_2)b\\
& &=(-1)^{\tilde {a} \tilde {c}}\Res_{z_{1}}\Res_{z_2}\frac{(1+z_{1})^{{\wt}c-\varepsilon}}
{z_{1}}\frac{(1+z_2)^{{\wt}a-1+\delta_{r}+\frac{r}{T}}}{z_2^{1+\delta_{r}}}Y(a,z_2)Y(c,z_{1})b\nonumber\\
& & \ +\Res_{z_{1}}\Res_{z_2}\Res_{z_{0}}
\frac{(1+z_{1})^{{\wt}c-\varepsilon}}{z_{1}}\frac{(1+z_2)^{{\wt}a-1+\delta_{r}+\frac{r}{T}}}{z_2^{1+\delta_{r}}}
z_2^{-1}\delta\left(\frac{z_{1}-z_{0}}{z_2}\right)Y(Y(c,z_{0})a,z_2)b\nonumber\\
& &=(-1)^{\tilde {c} \tilde {a}}\Res_{z_{1}}\Res_{z_2}\frac{(1+z_{1})^{{\wt}c-\varepsilon}}
{z_{1}}\frac{(1+z_2)^{{\wt}a-1+\delta_{r}+\frac{r}{T}}}{z_2^{1+\delta_{r}}}Y(a,z_2)Y(c,z_{1})b\nonumber\\
& &\ \ \ +\Res_{z_2}\Res_{z_{0}}
\frac{(1+z_2+z_{0})^{{\wt}c-\varepsilon}}{z_2+z_{0}}\frac{(1+z_2)^{{\wt}a-1+\delta_{r}+\frac{r}{T}}}{z_2^{1+\delta_{r}}}Y(Y(c,z_{0})a,z_2)b\nonumber\\
& &=(-1)^{\tilde {c} \tilde {a}}\Res_{z_2}\frac{(1+z_2)^{{\wt}a-1+\delta_{r}+\frac{r}{T}}}{z_2^{1+\delta_{r}}}Y(a,z_2)
\Res_{z_{1}}\frac{(1+z_{1})^{{\wt}c-\varepsilon}}
{z_{1}}Y(c,z_{1})b\nonumber\\
& &\ \ \ +\sum_{i,j=0}^{\infty}(-1)^{j}\left(\begin{array}{c}{\wt}c-\varepsilon\\i\end{array}\right)\Res_{z_2}
\frac{(1+z_2)^{{\wt}a-1+\delta_{r}+\frac{r}{T}+{\wt}c-\varepsilon -i}}
{z_2^{j+2+\delta_{r}}}Y(c_{i+j}a,z_2)b\nonumber\\
& &=(-1)^{\tilde {c} \tilde {a}}\Res_{z_2}\frac{(1+z_2)^{{\wt}a-1+\delta_{r}+\frac{r}{T}}}{z_2^{1+\delta_{r}}}Y(a,z_2)
\Res_{z_{1}}\frac{(1+z_{1})^{{\wt}c-\varepsilon}}
{z_{1}}Y(c,z_{1})b\nonumber\\
& &\ \ \ +\sum_{i,j=0}^{\infty}(-1)^{j}\left(\begin{array}{c}{\wt}c-\varepsilon\\i\end{array}\right)\Res_{z_2}
\frac{(1+z_2)^{{\wt}(c_{i+j}a)-1+\delta_{r}+\frac{r}{T}+j+1-\varepsilon}}
{z_2^{j+2+\delta_{r}}}Y(c_{i+j}a,z_2)b.\nonumber
\end{eqnarray*}
It is clear by the definition of $O_g(V)$ and 
Lemma \ref{l3.2} (i) that the  resulting vector lies  in $I.$
This shows that $I$ is an ideal of $V^{0*}.$  

Next we prove that $*_g$ is associative. We need to verify that
$(a*b)*c-a*(b*c)\in O_g(V^{0*})$ for $a,b,c \in V^{0*}.$ A straightforward
computation using the twisted Jacobi identity gives
 \begin{eqnarray*}
  & & (a*b)*c=\sum^{\wt a}_{i=0}(a_{i-1}b)*c \\
  &         & =\sum^{\wt a}_{i=0}\binom{\wt a}{i} \Res_w(Y(a_{i-1}b,w)\frac {(1+w)^{\wt (a_{i-1}b)}}{w}c)\nonumber\\
  &         & =\sum^{\wt a}_{i=0}\binom{\wt a}{i} \Res_w \Res_{z-w}(Y(Y(a,z-w)b,w)(z-w)^{i-1}\frac{(1+w)^{\wt a+\wt b-i}}{w} c)\\
  &         & = \Res_w \Res_{z-w}(Y(Y(a,z-w)b,w)\frac{(1+z)^{\wt a}(1+w)^{\wt b}}{w(z-w)}c)\\
  &         & =\Res_w \Res_z(Y(a,z)Y(b,w)\frac{(1+z)^{\wt a}(1+w)^{\wt b}}{w(z-w)}c)\\
  &         &-(-1)^{\tilde {a} \tilde {b}}\Res_w \Res_z(Y(b,w)Y(a,z)\frac{(1+z)^{\wt a}(1+w)^{\wt b}}{w(z-w)}c)\\
  &         & = \sum^{\infty}_{i=0}\Res_w \Res_z(Y(a,z)Y(b,w)(-1)^{i}z^{-1-i}w^i\frac{(1+z)^{\wt a}(1+w)^{\wt b}}{w}c)\\
  &            & -(-1)^{\tilde {a} \tilde {b}}\sum^{\infty}_{i=0}\Res_w \Res_z(Y(b,w)Y(a,z)(-1)^{i+1}z^{i}w^{-i-1}\frac{(1+z)^{\wt a}(1+w)^{\wt b}}{w}c)\\
  &         &\equiv \Res_z\Res_w(Y(a,z)Y(b,w)\frac{(1+z)^{\wt a}(1+w)^{\wt b}}{zw}c)\ \ mod \ \ O_g(V^{0*}) \\
 &          & \equiv a*(b*c) \ \ mod \ \ O_g(V^{0*})
\end{eqnarray*} 
Thus $A_g(V)\simeq \frac{V^{0*}}{O_g(V)\cap V^{0*}}$ is an associative algebra under $*_g.$ This finishes the proof of (i). The proofs of (ii) and (iii)
are immediate. \qed

\section{The Lie superalgebra $V[g]$}
\setcounter{equation}{0}
Let $V$\ \ be\ a\ vertex\ operator superalgebra with automorphism $g$ of order
$T_0.$ We can endow ${\C}[t^{\frac{1}{T_0}},t^{-\frac{1}{T_0}}]$ the structure of a vertex algebra with vertex operator
\begin{equation}\label{g4.1}
Y(f(t),z)g(t)=f(t+z)g(t)=\left(e^{z\frac{d}{dt}}f(t)\right)g(t).
\end{equation} 
(see \cite{B}). Then the tensor product
\begin{equation}\label{g4.2}
{\cal L}(V)={\C}[t^{\frac{1}{T_0}},t^{-\frac{1}{T_0}}]\otimes V.
\end{equation}
is a vertex superalgebra with vertex operator
\begin{equation}\label{g4.3}
Y(f(t)\otimes v,z)(g(t)\otimes u)=f(t+z)g(t)\otimes Y(v,z)u.
\end{equation} 
The $L(-1)$ operator of ${\cal L}(V)$ is given
by $D=\frac{d}{dt}\otimes 1+1\otimes L(-1).$ Extends 
$g$ to  an automorphism of vertex 
superalgebra in the following way:
\be{g4.3'}
g(t^{m}\otimes a)=\exp (\frac{-2 \pi im}{T_0})(t^{m}\otimes ga).
\end{equation}
Let ${\cal L}(V,g)$ be the $g$-invariants
which is a vertex sub-superalgebra of ${\cal L}(V).$ 
Clearly,
\be{g4.4}
{\cal L}(V,g)=\oplus_{r=0}^{T_{0}-1}t^{r/T_{0}}\C[t,t^{-1}]\otimes V^r.
\end{equation}

Following [B], we know that 
\be{g4.5}
V[g]={\cal L}(V,g)/D{\cal L}(V,g)
\end{equation}
is a  Lie superalgebra with bracket
\be{g4.6}
[u+D{\cal L}(V,g),v+D{\cal L}(V,g)]=u_{0}v+D{\cal L}(V,g).
\end{equation}

For short let $a(q)$ be  the image of $t^q\otimes a\in
{\cal L}(V,g)$ in $V[g].$ Then we have

\bl{l4.1} Let $a\in V^r,$ $v\in V^s$ and $m,n\in\Z.$ Then

(i) $[\omega(0),a(m+\frac{r}{T_{0}})]=-\left(m+\frac{r}{T_{0}}\right)a(m-1+\frac{r}{T_{0}}).$

(ii) $[a(m+\frac{r}{T_{0}}), b(n+\frac{s}{T_{0}})]=
\sum_{i=0}^{\infty}\binom{m+\frac{r}{T_0}}{i}a_ib(m+n+\frac{r+s}{T_{0}}-i).$

(iii) ${\1}(-1)$ lies in the center of $V[g].$
\el

For homogeneous $a\in V,$ we define 
\be{g4.7}
\deg(t^n\otimes a)=\wt a-n-1
\end{equation}
then ${\cal L}(V)$ is $\frac{1}{T}{\Z}$-graded. 
Since $D$ increases degree by 1, $D{\cal L}(V,g)$ is a graded
subspace of ${\cal L}(V,g)$ and $V[g]$ is naturally 
$\frac{1}{T}\Z$-graded:
$$V[g]=\oplus_{n\in\frac{1}{T}\Z}V[g]_n.$$ 
By Lemma \ref{l4.1}, $V[g]$ is a $\frac{1}{T}\Z$-graded Lie superalgebra 
with the  triangular decomposition
\be{g4.8}
V[g]=V[g]_{+}\oplus V[g]_0\oplus V[g]_{-}.
\end{equation} 
where 
$\displaystyle{V[g]_{\pm}=\sum_{0<n\in \frac{1}{T}{\Z}}V[g]_{\pm n}}.$

\begin{lem}\label{add} $V[g]_0$ is spanned by elements of the form 
$a(\wt a-1)$ for homogeneous $a\in V^{0*}.$ 
\end{lem}

\pf Let $a\in V.$ Then the degree $\wt a-n-1$ of  $a(n)$ is 0 
if and only if $a\in V_{\bar 0}^0$ and $n=\wt a-1$ or 
$a\in V_{\bar 1}^{T_0/2}$ and $n=\wt a-1.$ \qed

The bracket of $V[g]_0$ is given
by
\begin{eqnarray}\label{ec0}
[a({\wt}a-1),b({\wt}b-1)]=\sum_{j=0}^{\infty}\binom{{\wt}a-1}{j}a_jb({\wt}(a_{j}b)-1).
\end{eqnarray}

Set $o(a)=a(\wt a-1)$ for homogeneous $a\in V^{0*}$ and extend
linearly to all $a\in  V^{0*}.$ This gives a linear map 
\begin{eqnarray}
 V^{0*}&\to& V[g]_0,\nonumber\\
a &\mapsto& o(a).\label{g4.10}
\end{eqnarray}
As the  kernel of the map is $(L(-1)+L(0))V^{0*},$ 
we obtain an isomorphism of Lie superalgebras $V^{0*}/(L(-1)+L(0))V^{0*}\cong V[g]_0.$ The  bracket on the quotient of $V^{0*}$  is given by
$$[a,b]=\sum_{j\geq 0}\binom{{\wt}a-1}{j}a_{j}b.$$

\bl{l4.2} Let $A_g(V)_{Lie}$ be the Lie superalgebra of the associative algebra
$A_g(V)$ introduced in section 3 such that  
$[u,v]=u*_gv-(-1)^{\tilde {u}\tilde{v}}v*_gu.$ Then the map 
$o(a)\mapsto a+O_g(V)$ is an onto  Lie superalgebra homomorphism from
$V[g]_0$ to $A_g(V)_{Lie}.$
\el

\pf Recall that $I=O_g(V)\cap V^{0*}.$ So we have a surjective linear map
$$V[g]_0\cong V^{0*}/(L(-1)+L(0))V^{0*}\to  V^{0*}/I \simeq A_g(V),$$
\be{g4.12}
o(a)\to a+(L(-1)+L(0))V^{0*}\to a+ I.
\end{equation}
The Lie homomorphism follows from
 $$[o(a),o(b)]=\sum^{\infty}_{j=0}\binom{\wt a -1}{j}o(a_jb).$$
and
\begin{eqnarray*}
 [a+O_g(V),b+O_g(V)]&\equiv& a*_gb-(-1)^{\tilde{a} \tilde{b}}b*_ga \\
                      &\equiv&\sum^{\infty}_{j=0}\binom{\wt a-1}{j}a_jb\\
                      &\equiv& \Res_z(1+z)^{\wt a-1}Y(a,z)b \ \ mod\ \ O_g(V^{0*})\\
         &\equiv&\sum^{\infty}_{i=0}\binom{\wt a -1}{i}a_ib \ \ mod \ \ O_g(V^{0*}).
\end{eqnarray*}
\qed

\section{The functor $\O$}
\setcounter{equation}{0}
The main purpose in this section is to construct a covariant functor $\Omega$ from the category of weak $g$-twisted $V$-modules to the category
of $A_g(V)$-modules (cf. Theorem \ref{t5.3}). Let $M$ be a weak
$g$-twisted $V$-module. We define the space of ``lowest weight vectors'' to be 
$$\Omega(M)=\{w\in M|u_{\wt u+n}w=0, u\in V,n\geq 0\}.$$ 
The main result 
in this section says that $\Omega(M)$ is an $A_g(V)$-module.  
Moreover if $f:M\to N$ is a morphism in weak $g$-twisted $V$-modules,
the restriction $\Omega(f)$ of $f$ to $\O(M)$ is an $A_g(V)$-module morphism.

Note that if $M$ is a weak $g$-twisted $V$-module then 
$M$ becomes a $V[g]$-module such that $a(m)$ acts as $a_m.$ 
Moreover, $M$ is an admissible $g$-twisted
$V$-module if and only if $M$ is a $\frac{1}{T}\Z_+$-graded module for
the graded Lie superalgebra $V[g]$ (cf. Lemmas 5.1 and 5.2 of \cite{DLM2}).

\bt{t5.3} Let  $M$ be  a weak $g$-twisted $V$-module. Then the map $a\mapsto o(a)$  for homogeneous $a\in V^{0*}$ gives a representation of the associative algebra $A_g(V)$ on $\Omega(M)$.
\et

\pf We first show that $\Omega(M)$ is invariant under the action of
$o(a)$ for $a\in V.$ Let $b\in V,$ $w\in M$ and $n\geq \wt b.$  By (\ref{g2.17}),
$$b_no(a)v=b_na_{\wt a-1}v=a_{\wt a-1}b_nv+\sum_{i\geq 0}\binom{n}{i}
(b_ia)_{\wt a+n-1-i}v=0$$
as  $\wt a+n-1-i$ is greater than or equal to  $\wt (b_ia).$ 

Next we show that the action is well defined. For $a\in V^{r*}$
and $0<r<T,$  $o(a)=0$ by Lemma \ref{add}.
It remains to check that $o(a)=0$ on $\Omega(M)$ for $a\in I.$

Note that $a(\wt a-1+n)\Omega(M)=0$ for any $n>0$ and $a\in V^{0*}$.
There are two cases: $a\in O(V^{0*})$ or  
$$a=\Res_{z}\frac{(1+z)^{{\wt}c-1+\frac{r}{T}}}{z}Y(u,z)v$$
with $u\in V^{r*}, v\in V^{(T-r)*}, 0\le r\le T-1$.

If $a\in O(V^{0*})$ then there exist $u,v\in V^{0*}$ such that
$$a=Res_z\frac{(1+z)^{\wt u}}{z^{2}}Y(u,z)v.$$
The argument in the Proof of Theorem 2.1.2 in \cite{Z} with suitable
modification gives
$$o(u*v)=o(u)o(v).$$
Note that $o(L(-1)u+L(0)u)=0$ and $(L(-1)u+L(0)u)*v=u\circ v.$ We immediately
have $o(u\circ v)=0$ on $\Omega(M).$

If 
$$a=\Res_{z}\frac{(1+z)^{{\wt}c-1+\frac{r}{T}}}{z}Y(u,z)v,$$
we can use Lemma \ref{ladd}. Since $z^{\wt u-1+\frac{r}{T}}Y_M(u,z)w$
involves only nonnegative integer powers of $z$ for $w\in \Omega(M),$
we have  
\begin{equation}\label{eqn}
(z_{0}+z_{2})^{\wt u-1+\frac{r}{T}}Y_{M}(u,z_{0}+z_{2})Y_{M}(v,z_{2})w
=(z_{2}+z_{0})^{\wt u-1+\frac{r}{T}}Y_M(Y(u,z_0)v,z_2)w.
\end{equation}
Applying $\Res_{z_{0}}\Res_{z_{2}}z_{0}^{-1}z_{2}^{{\wt}v-\frac{r}{T}}$ to (\ref{eqn}) yields
\begin{eqnarray} 
& &0=\Res_{z_{0}}\Res_{z_{2}}z_{0}^{-1}z_{2}^{{\wt}v-\frac{r}{T}}
(z_{2}+z_{0})^{{\wt}u-1+\frac{r}{T}}Y_{M}(Y(u,z_{0})v,z_{2})w\nonumber\\
& &\ \ \ \ =\sum_{i=0}^{\infty}\left(\begin{array}{c}{\wt}u-1+\frac{r}{T}\\i\end{array}\right)
\Res_{z_{2}}z_{2}^{{\wt}u+{\wt}v-i-1}Y_{M}(u_{i-1}v,z_{2})w\nonumber\\
& &\ \ \ \ =\sum_{i=0}^{\infty}\left(\begin{array}{c}{\wt}u-1+\frac{r}{T}\\i\end{array}\right)o(u_{i-1}v)w
\nonumber\\
& &\ \ \ \ =o\left(\Res_{z}\frac{(1+z)^{{\wt}u-1+\frac{r}{T}}}{z}Y_{M}(u,z)v\right)w\nonumber\\
& &\ \ \ \ =o(a)w
\end{eqnarray}
as required. 
\qed

If $M$ is a nonzero   admissible $g$-twisted $V$-modules we may and
do assume that $M(0)$
is nonzero with suitable degree shift. With these
conventions we have
 
\bp{ll2.9} Let  $M$ be a simple  admissible $g$-twisted $V$-module. 
Then the following hold

(i) $\Omega(M)=M(0).$ 

(ii) $\O(M)$ is a simple $A_{g}(V)$-module.
\ep

\pf The proof is the  same as in [DLM2].

\section{Generalized Verma modules and the functor $L$}
\setcounter{equation}{0}
In this section
we focus on how to  construct admissible $g$-twisted $V$-modules
from a given $A_g(V)$-module $U.$  We use the same trick which was used in 
[DLM2] to do this. We will define two $g$-twisted admissible $V$-modules
$\bar M(U)$ and $L(U).$  The $\bar M(U)$ is  the universal
admissible $g$-twisted $V$-module such that $\bar M(U)(0)=U$
and $L(U)$ is smallest admissible $g$-twisted $V$-module 
whose $L(U)(0)=U.$ Just as in the classical highest weight module theory,
$L(U)$ is the unique irreducible quotient of $\bar M(U)$ if $U$ is simple.

We start with an  $A_{g}(V)$-module $U.$
Then $U$ is automatically a module for $A_g(V)_{Lie}.$
By  Lemma \ref{l4.2} $U$ is lifted to a module for the Lie superalgebra
$V[g]_0.$ Let $V[g]_-$ act trivially on $U$ and extend 
$U$ to a $P=V[g]_{-}\oplus V[g]_{0}$-module.
 Consider the induced module 
\be{g6.1}
M(U)=\Ind_{P}^{V[g]}(U)=U(V[g])\otimes_{U(P)} U 
\end{equation}
which is $\frac{1}{T}\Z_+$-graded module for $V[g]$ by giving 
$U$ degree 0. Then  $M(U)(n)=U(V[g]_+)_nU$ by PBW theorem and in
particular $M(U)(0)=U.$ 

For $v\in V$ we set 
\begin{equation}\label{g6.2}
Y_{M(U)}(v,z)=\sum_{m\in\frac{1}{T_0}\Z}v(m)z^{-m-1}
\end{equation}
Then $Y_{M(U)}(v,z)$ satisfies condition (\ref{1/2})-(\ref{vacuum}).
By Lemma \ref{l4.1} (ii), the identity (\ref{g2.17})
holds. But this is not good enough to establish the twisted Jacobi identity 
for the action (\ref{g6.2}) on $M(U).$ 

Let $W$ be the subspace of $M(U)$ spanned linearly by the 
coefficients of
\begin{eqnarray}\label{g6.3}
(z_{0}+z_{2})^{{\wt}a-1+\delta_r+\frac{r}{T}}Y(a,z_{0}+z_{2})Y(b,z_{2})u-(z_{2}+z_{0})^{{\wt}a-1+\delta_r+\frac{r}{T}}
Y(Y(a,z_{0})b,z_{2})u
\end{eqnarray}
for any homogeneous $a\in V^{r*},b\in V,$ $u\in U$.
We set
\be{g6.4}
\bar M(U)=M(U)/U(V[g])W.
\end{equation}

\bp{p6.1}  Let $M$ be a $V[g]$-module such that there is a 
 subspace $U$ of $M$ satisfying the
following conditions:

(i) $M=U(V[g])U;$

(ii) For any $a\in V^{r*}$ and $u\in U$ there is 
$k\in \wt\,a+\Z_+$ 
such that
\begin{eqnarray}\label{d1}
(z_{0}+z_{2})^{k+\frac{r}{T}}Y(a,z_{0}+z_{2})Y(b,z_{2})u=
(z_{0}+z_{2})^{k+\frac{r}{T}}Y(Y(a,z_{0})b,z_{2})u   
\end{eqnarray}
for any $b\in V$. Then $M$ is a weak $V$-module.
\ep

\pf We only need to prove the twisted Jacobi identity , which is equivalent
to commutator relation (\ref{g2.16}) and  the associativity
(\ref{ea2.15}). But the commutator formula is built in already as
$M$ is a $V[g]$-module.  By Lemma \ref{ladd}, the assumption
(ii) can be reformulated as follows:

(ii') For any $a\in V^{r}$ and $u\in U$ there is 
$k\in \Z_+$ 
such that
\begin{eqnarray}\label{d2}
(z_{0}+z_{2})^{k+\frac{r}{T_0}}Y(a,z_{0}+z_{2})Y(b,z_{2})u=
(z_{0}+z_{2})^{k+\frac{r}{T_0}}Y(Y(a,z_{0})b,z_{2})u   
\end{eqnarray}

Since $M$ is a $V[g]$-module generated by $U$ it is enough to prove that
if $u$ satisfies (ii') then $c_nu$ also satisfies (ii') for
$c\in V$ and $n\in \frac{1}{T_0}\Z.$  

Let $k_{1}$ be a positive integer such that $c_{i}a=0$ 
for $i\ge k_{1}$. Let $k_2$ be a positive integer such that
\begin{eqnarray}
& &(z_{0}+z_{2})^{k_{2}+\frac{r+s}{T_0}}Y(c_{i}a,z_{0}+z_{2})Y(b,z_{2})u
=(z_{2}+z_{0})^{k_{2}+\frac{r+s}{T_0}}Y(Y(c_{i}a,z_{0})b,z_{2})u,\label{3.18}\\
& &(z_{0}+z_{2})^{k_{2}+\frac{r+s}{T_0}}Y(a,z_{0}+z_{2})Y(c_{i}b,z_{2})u
=(z_{2}+z_{0})^{k_{2}+\frac{r+s}{T_0}}Y(Y(a,z_{0})c_{i}b,z_{2})u\hspace{1cm}\label{3.19}
\end{eqnarray}
for any nonnegative integer $i$ where we have assumed that 
$c\in V^s.$ 

Let $k$ be a positive integer such that
$k+\frac{r}{T_0}+n-k_{1}>k_{2}+\frac{r+s}{T_0}.$ Using (\ref{3.18})  and
(\ref{3.19}) and the bracket relation (ii) in Lemma \ref{l4.1}
\begin{eqnarray}\label{3.6'}
[a_m,Y(b,z_{2})]
=\sum_{i=0}^{\infty}\binom{m}{i}z_2^{m-i}Y(a_ib,z_2)
\end{eqnarray}
we have
\begin{eqnarray*}
& &\ \ \ \ \ (z_{0}+z_{2})^{k+\frac{r}{T_0}}Y(a,z_{0}+z_{2})Y(b,z_{2})c_{n}u\\
& &=(-1)^{\tilde {a}\tilde {c}}(-1)^{\tilde {b}\tilde {c}}(z_{0}+z_{2})^{k+\frac{r}{T_0}}c_{n}Y(a,z_{0}+z_{2})Y(b,z_{2})u\\
& &\ \ \ \ -(-1)^{\tilde {a}\tilde {c}}(-1)^{\tilde {b}\tilde {c}} \sum_{i=0}^{\infty}\binom{n}{i}
(z_{0}+z_{2})^{k+\frac{r}{T_0}+n-i}Y(c_{i}a,z_{0}+z_{2})Y(b,z_{2})u\\
& &\ \ \ \  -(-1)^{\tilde {b}\tilde {c}}\sum_{i=0}^{\infty}\binom{n}{i}
z_{2}^{n-i}(z_{0}+z_{2})^{k+\frac{r}{T_0}}Y(a,z_{0}+z_{2})Y(c_{i}b,z_{2})u\\
& &=(-1)^{\tilde {a} \tilde {c}}(-1)^{\tilde {b}\tilde {c}}(z_{0}+z_{2})^{k+\frac{r}{T_0}}c_{n}Y(Y(a,z_{0})b,z_{2})u\\
& &\ \ \ \ -(-1)^{\tilde {a} \tilde {c}}(-1)^{\tilde {b} \tilde {c}}\sum_{i=0}^{\infty}\binom{n}{i}
(z_{2}+z_{0})^{k+\frac{r}{T_0}+n-i}Y(Y(c_{i}a,z_{0})b,z_{2})u\\
& &\ \ \ \ -(-1)^{\tilde {b}\tilde {c}}\sum_{i=0}^{\infty}\binom{n}{i}
z_{2}^{n-i}(z_{2}+z_{0})^{k+\frac{r}{T_0}}Y(Y(a,z_{0})c_{i}b,z_{2})u\\
& &=(-1)^{\tilde {a} \tilde {c}}(-1)^{\tilde {b}\tilde {c}}(z_{0}+z_{2})^{k+\frac{r}{T_0}}c_{n}Y(Y(a,z_{0})b,z_{2})u\\
& &\ \ \ \ -(-1)^{\tilde {a} \tilde {c}}(-1)^{\tilde {b}\tilde {c}}\sum_{i=0}^{\infty}\binom{n}{i}
(z_{2}+z_{0})^{k+\frac{r}{T_0}+n-i}Y(Y(c_{i}a,z_{0})b,z_{2})u\\
& &\ \ \ \ -(-1)^{\tilde {a} \tilde {c}}(-1)^{\tilde {b}\tilde {c}}\sum_{i=0}^{\infty}\binom{n}{i}
z_{2}^{n-i}(z_{2}+z_{0})^{k+\frac{r}{T_0}}Y(c_{i}Y(a,z_{0})b,z_{2})u\\
& &\ \ \ \ +(-1)^{\tilde {a} \tilde {c}}(-1)^{\tilde {b}\tilde {c}}\sum_{i=0}^{\infty}\sum_{j=0}^{\infty}\binom{n}{j}\binom{j}{i}
z_{2}^{n-i}(z_{2}+z_{0})^{k+\frac{r}{T_0}}z_{0}^{i-j}Y(Y(c_{j}a,z_{0})b,z_{2})u\\
& &=(-1)^{\tilde{b}\tilde{c}}(-1)^{\tilde{a}\tilde{c}}(z_{0}+z_{2})^{k+\frac{r}{T_0}}c_{n}Y(Y(a,z_{0})b,z_{2})u\\
& &\ \ \ \ -(-1)^{\tilde {a} \tilde {c}}(-1)^{\tilde {b}\tilde {c}}\sum_{i=0}^{\infty}\binom{n}{i}
(z_{2}+z_{0})^{k+\frac{r}{T_0}+n-i}Y(Y(c_{i}a,z_{0})b,z_{2})u\\
& &\ \ \ \ -(-1)^{\tilde {a} \tilde {c}}(-1)^{\tilde {b}\tilde {c}}\sum_{i=0}^{\infty}\binom{n}{i}
z_{2}^{n-i}(z_{2}+z_{0})^{k+\frac{r}{T_0}}Y(c_{i}Y(a,z_{0})b,z_{2})u\\
& &\ \ \ \ +(-1)^{\tilde {a} \tilde {c}}(-1)^{\tilde {b}\tilde {c}}\sum_{j=0}^{\infty}\sum_{i=j}^{\infty}\binom{n}{j}\binom{n-j}{i-j}
z_{2}^{n-i}(z_{2}+z_{0})^{k+\frac{r}{T_0}}z_{0}^{i-j}
Y(Y(c_{j}a,z_{0})b,z_{2})u \nonumber\\
& &=(-1)^{\tilde {a} \tilde {c}}(-1)^{\tilde {b}\tilde {c}}(z_{0}+z_{2})^{k+\frac{r}{T_0}}c_{n}Y(Y(a,z_{0})b,z_{2})u\\
& &\ \ \ \ -(-1)^{\tilde {a} \tilde {c}}(-1)^{\tilde {b}\tilde {c}}\sum_{i=0}^{\infty}\binom{n}{i}
z_{2}^{n-i}(z_{2}+z_{0})^{k+\frac{r}{T_0}}Y(c_{i}Y(a,z_{0})b,z_{2})u\\
& &=(-1)^{\tilde{b}\tilde{c}}(-1)^{\tilde{a}\tilde{c}}(z_{0}+z_{2})^{k+\frac{r}{T_0}}c_{n}Y(Y(a,z_{0})b,z_{2})u\\
& &\ \ \  \ -(-1)^{\tilde {a} \tilde {c}}(-1)^{\tilde {b}\tilde {c}}(z_{2}+z_{0})^{k+\frac{r}{T_0}}[c_n, Y(Y(a,z_{0})b,z_{2})]u\\
& &=(-1)^{\tilde{b}\tilde{c}}(-1)^{\tilde{a}\tilde{c}}(-1)^{(\tilde{a}+\tilde{b})\tilde{c}}(z_{2}+z_{0})^{k+\frac{r}{T_0}}Y(Y(a,z_{0})b,z_{2})c_{n}u,\\
& &=(z_{2}+z_{0})^{k+\frac{r}{T_0}}Y(Y(a,z_{0})b,z_{2})c_{n}u.
\end{eqnarray*}
The proof is complete. \qed

Applying Proposition \ref{p6.1} to $\bar M(U)$  gives 
the  following main result of this section.
\bt{t6.1} $\bar M(U)$ is an admissible $g$-twisted $V$-module with
$\bar M(U)(0)=U$ and with the 
following universal property: for any weak $g$-twisted $V$-module $M$
and any $A_g(V)$-morphism $\phi: U\to \O(M),$ there is a unique morphism 
$\bar\phi: \bar M(U)\to M$ of weak $g$-twisted $V$-modules which 
extends $\phi.$ 
\et

As in \cite{DLM2} we also have
\bt{t6.3} $M(U)$ has a unique maximal graded $V[g]$-submodule $J$ 
with the property that $J\cap U=0.$ Then $L(U)=M(U)/J$ is an admissible
$g$-twisted $V$-module satisfying $\O(L(U))\cong U.$

$L$ defines a functor from the category of 
$A_{g}(V)$-modules to
the category of admissible $g$-twisted $V$-modules such that $\O\circ L$
is naturally equivalent to the identity.
\et

We have a pair of functors $\O,L$ between the $A_g(V)$-module category
and admissible $g$-twisted $V$-module category.
Although $\O\circ L$ is equivalent to the identity, $L\circ\O$ is not
equivalent to the identity in general. 

The following result is an immediate consequence of Theorem \ref{t6.3}.

\bl{l7.1} Suppose that $U$ is a simple $A_{g}(V)$-module. 
Then $L(U)$ is a simple admissible $g$-twisted $V$-module.
\el

Using Lemma \ref{l7.1}, Proposition 
\ref{ll2.9} (ii), Theorems \ref{t6.1} and \ref{t6.3} gives: 
\bt{t7.2} $L$ and $\O$ are equivalent when restricted to the full 
subcategories of completely reducible $A_g(V)$-modules and 
completely reducible admissible $g$-twisted $V$-modules respectively.
In particular, $L$ and $\Omega$ induces mutually
inverse bijections on the isomorphism classes of simple objects
in the category of $A_{g}(V)$-modules and admissible $g$-twisted $V$-modules
respectively.
\end{thm}

We now apply the obtained results to $g$-rational vertex operator 
superalgebras to obtain:

\bt{t8.1} Suppose that $V$ is a $g$-rational vertex operator superalgebra. Then the 
following hold:

(a) $A_g(V)$ is a finite-dimensional, semi-simple associative algebra
(possibly 0).

(b) $V$ has only finitely many isomorphism classes of simple admissible
 $g$-twisted modules.

(c) Every simple admissible $g$-twisted $V$-module is an ordinary $g$-twisted
$V$-module.

(d) $V$ is $g^{-1}$-rational.

(e) The functors $L,\O$ are mutually inverse categorical equivalences
between the category of $A_g(V)$-modules and the category of admissible $g$-twisted $V$-modules.

(f) The functors $L,\O$ induce mutually inverse categorical equivalences
between the category of finite-dimensional 
$A_g(V)$-modules and the category of ordinary $g$-twisted $V$-modules.
\et

The proof is the same as that of Theorem 8.1 in \cite{DLM2}.

\section{Examples}
\setcounter{equation}{0}

In this section we discuss the well known vertex operator superalgebras
constructed from the free fermions and their twisted modules. In 
particular we compute the algebra $A_g(V)$ and classify the irreducible
twisted modules using $A_g(V).$ The classification results
have  been obtained previously in \cite{Li2} with a different 
approach. 
 
Let $H=\sum_{i=1}^{l}{\Bbb C}a_{i}$ be a complex vector space equipped
with a nondegenerate symmetric bilinear form $(,)$ such that  $\{a_i|
i=1,2,...l\}$ form an orthonormal basis.
 Let $A(H,{\Bbb Z}+\frac{1}{2})$ be the
associative algebra generated by $\{a(n)|a\in H,n\in {\Bbb
Z}+\frac{1}{2}\}$ subject to the relation
$$[a(n),b(m)]_{+}=(a,b)\delta_{{m+n},0}.$$ 
Let $A^{+}(H,{\Bbb
Z} +\frac{1}{2})$ be the subalgebra generated by $\{a(n)|a\in H,n\in
{\Bbb Z}+\frac{1}{2},n>0 \},$ and make ${\Bbb C}$ a $1$-dimensional
$A^{+}(H,{\Bbb Z} +\frac{1}{2})$-module so that $a_i(n)1=0$
for $n>0$.

Consider the induced module
\begin{eqnarray*}
& &V(H,{\Bbb Z}+\frac{1}{2})=A(H,{\Bbb Z}+\frac{1}{2})\otimes_{A^{+}(H,{\Bbb Z}+\frac{1}{2})}\C\\
& & \ \ \ \ \ \ \ \ \ \  \ \ \ \ \  \  \ \cong \Lambda [a_i(-n)|n >0,n\in {\Bbb Z}+\frac{1}{2} ,i=1,2,...l]\ ({{\rm linearly}}).
 \end{eqnarray*}
The action of $a_i(n)$ is given by $\frac{\partial}{\partial a_i(-n)}$
if $n$ is positive and by multiplication by $a_i(n)$ if $n$ is negative.

The $V(H,{\Bbb Z}+\frac{1}{2})$ is naturally graded by $\frac{1}{2}\Z$
so that 
$$V(H,{\Bbb Z}+\frac{1}{2})_n= \<a_{i_1}(-n_1)a_{i_2}(-n_2)\cdots a_{i_k}(-n_k)|n_1+n_2+\cdots n_k=n \>$$
 
Let $b_1,...,b_k\in H$ and $n_1,...,n_k\in\frac{1}{2}\Z$ we define a normal
ordering
$$:b_1(n_1)\cdots b_k(n_k):=(-1)^{|\sigma|} b_{i_1}(n_{i_1})\cdots b_{i_k}(n_{i_k})$$
such that $n_{i_1}\leq \cdots \leq n_{i_k}$ where $\sigma$ is the 
permutation of $\{1,...,k\}$ by sending $j$ to $i_j.$  
For $a\in H$ set 
$$Y(a(-1/2),z)=\sum_{n\in\frac{1}{2}+\Z}a(n)z^{-n-1/2}.$$
Let $v=b_1(-n_1-\frac{1}{2})\cdots b_k(-n_k-\frac{1}{2})$ be
a general vector in $V(H, {\Bbb Z}+\frac{1}{2})$  where $n_i$ are nonnegative
integers. We set
$$Y(v,z)=:(\partial_{n_1}b_1(z))\cdots (\partial_{n_k}b_k(z)):$$
where $\partial_n=\frac{1}{n!}(\frac{d}{dz})^n.$ Then we have
a linear map:
\begin{equation}
\begin{array}{ccc}
V(H,{\Bbb Z}+\frac{1}{2})&\to& (\mbox{End}\,V(H,{\Bbb Z}+\frac{1}{2}))
[[z,z^{-1}]]\hspace*{3.6 cm} \\
v&\mapsto& Y(v,z)=\displaystyle{\sum_{n\in{\Bbb Z}}v_nz^{-n-1}\ \ \ (v_n\in
\mbox{End}\,V(H,{\Bbb Z}+\frac{1}{2}))}.
\end{array}
\end{equation}

Set $\1=1$ and  $\omega=\frac{1}{2}\sum_{i=1}^{l}a_i(-\frac{3}{2})a_i(-\frac{1}{2}).$ The following result is well known (cf. \cite{FFR}, 
 \cite{KW} and \cite{Li1}).
\begin{thm} $(V(H,{\Bbb Z}+\frac{1}{2}), Y, \1,\omega)$ is a holomorphic
vertex operator superalgebra generated by $a_i(-\frac{1}{2})$ for $i=1,...,l.$
\end{thm}

We have already mentioned in Section 2 that any vertex operator
superalgebra has a canonical automorphism $\sigma$ such that
$\sigma=1$ on $V_{\bar 0}$ and $\sigma=-1$ on $V_{\bar 1}.$ Note that
$ V(H,{\Bbb Z}+\frac{1}{2})_{\bar 1}\ne 0.$ So $\sigma $ is an order
2 automorphism of  $V(H,{\Bbb Z}+\frac{1}{2})_{\bar 1}.$ We next
discuss the $\sigma$-twisted  $V(H,{\Bbb Z}+\frac{1}{2})$-modules.

First we discuss the case  when $l=2k$ ($k$ is a positive integer) is even. 
The $H$ can be written as 
$$H=\sum_{i=1}^{k}{\Bbb C} b_i+\sum_{i=1}^{k}{\Bbb C}b_i^* $$ with
$(b_i,b_j)=(b_i^*,b_j^*)=0, (b_i,b_j^*)=\delta_{i,j}.$ Let $A(H,{\Bbb
Z})$ be the associative algebra generated by $\{b(n)|b\in H, n\in
{\Bbb Z}\} $ subject to the relation
$$[a(n),b(n)]_+=(a,b)\delta _{m+n,0}$$ Let $A(H,{\Bbb Z})^+$ be the
subalgebra generated by $\{b_i(n),b_i^{*}(m)|n\geq 0, m >0,
i=1,\cdots, k\},$ and make ${\Bbb C}$ a $1$-dimensional $A(H,{\Bbb
Z})^+$-module so that $b_i(n)1=0$ for $n>0 $ and $b_i^*(n)1=0$ for
$n\geq 0$ $i=1,\cdots, k.$

As before we consider the induced module 
\begin{eqnarray}
 & & V(H,{\Bbb Z})=A(H,{\Bbb Z})\otimes_{A(H,{\Bbb Z})^+} {\Bbb C}\\
 & &\ \ \ \ \ \ \ \ \ \ \  \ \ \ \cong \Lambda [b_i(-n),b_i^{*}(-m)|n,m\in {\Bbb Z},n>0, m\geq0].
\end{eqnarray} 
Then $b_i(n) $ acts as $\frac{\partial}{\partial b_i(-n)^*}$ if $n$ is
nonnegative and multiplication  by $b_i(n) $
if $n$ is negative. Similarly, $b_i(n)^* $ acts 
as  $\frac{\partial}{\partial b_i(-n)}$ if $n$ is
positive and multiplication by  $b_i(n)^* $ if $n$ is nonpositive.

Thanks to Proposition 4.3 in [Li2], $V(H,{\Bbb Z})$ is an irreducible
 $\sigma$-twisted $V(H,{\Bbb Z}+\frac{1}{2})$-module such that
\begin {eqnarray*}
Y_{V(H,{\Bbb Z})}(u(-\frac{1}{2}),z)=u(z)=\sum_{n\in {\Bbb Z}}u(n)z^{-n-1/2}
\end{eqnarray*}
for $u\in H.$

\begin {prop}\label{p7.1}
If $\dim H=l=2k$ is even, then $A_{\sigma}( V(H,{\Bbb Z}+\frac{1}{2})$
is isomorphic to the matrix algebra $M_{2^k\times 2^k}(\C)$ and 
$V(H,{\Bbb Z})$ is the unique irreducible
$\sigma$-twisted $ V(H,{\Bbb Z}+\frac{1}{2})$-module up to isomorphism.   
\end{prop}

\pf By theorem \ref{t5.3}, it is enough to show that 
$A_{\sigma}( V(H,{\Bbb Z}+\frac{1}{2})$ is isomorphic to the matrix algebra 
$M_{2^k\times 2^k}(\C).$

Since $g=\sigma,$ the decomposition (\ref{g2.9}) becomes 
$V=V^{0*}.$ By lemma \ref{l3.2}(i), 
\begin{eqnarray*}
\Res_z\frac{(1+z)^{\frac{1}{2}}}{z^{2+m}}a_i(z)v=\sum_{s\geq 0}c_sa_i(-m+s-\frac{3}{2})v
\end{eqnarray*}
lies in $O_{\sigma}(V(H,\Z+\frac{1}{2}))$ for $m\geq 0$ where
$c_s=\binom{\frac{1}{2}}{s}.$
This implies that  $$a_i(-m-\frac{3}{2})v \equiv \sum_{s=1}^{\infty}c_sa_i(-m-\frac{3}{2}+s)v \ \ \ mod \ \ \ O_{\sigma}(V(H,{\Bbb Z}+\frac{1}{2})).$$
Thus $A_{\sigma}(V(H,{\Bbb Z}+\frac{1}{2}))$ is spanned by
$b_1(-1/2)^{s_1}\cdots b_k(-1/2)^{s_k}b_1^*(-1/2)^{t_1}\cdots b_k^*(-1/2)^{t_k}$ with $s_i,t_i=0,1.$ As a result,  
$\dim A_{\sigma}(V(H,{\Bbb Z}+\frac{1}{2}))\leq 2^{2k}.$ 

On the other hand, $V(H,{\Bbb Z})$ is an irreducible $\sigma$-twisted
$V(H,{\Bbb Z}+\frac{1}{2})$-module. By Theorem \ref {t5.3}, $
\Omega(V(H,{\Bbb Z}))$ is a simple $A_{\sigma}(V(H,{\Bbb
Z}+\frac{1}{2}))$-module and $\dim\Omega(V(H,{\Bbb Z}))=2^k.$ Then $\dim
A_{\sigma}(V(H,{\Bbb Z}+\frac{1}{2}) \geq \dim \Omega(V(H,{\Bbb
Z}))=2^{2k}.$ This forces  $\dim A_{\sigma}(V(H,{\Bbb
Z}+\frac{1}{2}))=2^{2k}$ and $A_{\sigma}(V(H,{\Bbb Z}+\frac{1}{2}))\cong
M_{2^k\times 2^k}({\Bbb C}).$ 
\qed

We now deal with the case $\dim H=2k+1$ for some nonnegative integer $k.$
Then  $H$ can be decomposed into:
$$H=\sum_{i=1}^{k}{\Bbb C} b_i+\sum_{i=1}^{k}{\Bbb C}b_i^*+{\Bbb C}e$$ 
with $(b_i,b_j)=(b_i^*,b_j^*)=0,  (b_i,b_j^*)=\delta_{i,j}, (e,b_i)=(e,b_i^*)=0,(e,e)=2.$  

Let $A(H,{\Bbb Z})$ be the associative algebra generated same as above,
and $A(H,{\Bbb Z})^+$ be the subalgebra generated by $\{b_i(n),b_i^*(m),e(n)|m,n\in \Z,m>0,n\geq 0, i=1,\cdots, k\}$ and make ${\Bbb C}$  a $1$-dimensional $A(H,{\Bbb Z})^+$-module so that $b_i(n)1=0$ for $n\geq0 $ and $b_i^*(m)1=e(m)1=0$ for $m>0,$ $i=1,\cdots, k.$
Set  
\begin{eqnarray*}
 & & V(H,{\Bbb Z})=A(H,{\Bbb Z})\otimes_{A(H,{\Bbb Z})^+} {\Bbb C}\\
 & &\ \ \ \ \ \ \ \ \ \ \  \ \ \ \cong \Lambda [b_i(-n),b_i^{*}(-m),e(-m)|n,m\in {\Bbb Z},n>0, m\geq0]
\end{eqnarray*} 
and let
$$ W(H,{\Bbb Z})=\Lambda [b_i(-n),b_i^{*}(-m),e(-n)|n,m\in {\Bbb
    Z},n>0, m\geq0]=W(H,{\Bbb Z})^{even}\oplus W(H,\Z)^{odd}$$
be the decomposition into the even and old parity subspaces. 
Also define
$$V_{\pm}(H,\Z)=(1\pm e(0))W(H,{\Bbb Z})^{even}\oplus (1\mp e(0))W(H,{\Bbb Z})^{odd}.$$ 
Then 
$$V(H,\Z)=V_+(H,\Z)\oplus V_-(H,\Z)$$
and $V_{\pm}(H,\Z)$ are irreducible $A(H,{\Bbb Z})$-modules. The actions
of $b_i(n),b_i^*(n)$ are the same as before. The $e(n)$ acts as 
$2\frac{\partial}{\partial e(-n)}$ if $n>0$ and as multiplication
by $e(n)$ if $n\leq 0.$ 

Again by Proposition 4.3 in \cite{Li2},  $V_{\pm}(H,\Z)$ 
are irreducible $\sigma$-twisted  modules for $V(H,{\Bbb Z+1/2})$ so
that 
\begin {eqnarray*}
Y_{V(H,{\Bbb Z})}(u(-\frac{1}{2}),z)=u(z)=\sum_{n\in {\Bbb Z}}u(n)z^{-n-1/2}
\end{eqnarray*}
for $u\in H.$

\begin {prop}\label{p7.2}
If $\dim H=2k+1$ is odd, then  $ A_{\sigma}(V(H,{\Bbb Z}+\frac{1}{2}))$ is 
isomorphic to $M_{2^k\times 2^k}(\C)\oplus M_{2^k\times 2^k}(\C)$
and $V(H,{\Bbb Z}+\frac{1}{2})$  has exactly
two irreducible $\sigma$-twisted  modules $V_{\pm}(H,\Z)$ 
up to isomorphism.    
\end{prop}

\pf The proof is similar to that of Proposition \ref{p7.1}.\qed

Note that the automorphism $\sigma$ of $V(H,{\Bbb Z}+\frac{1}{2})$
is a lifting of $-1$-isometry of $H.$ It is natural to
study the twisted modules for an arbitrary isometry of $H.$ 

Let $H$ be a complex vector space as before with nondegenerate
symmetric bilinear form $(,)$ and $\tau$ an isometry of $H$ of
order $N_0<\infty.$ Then we can extend $\tau$ to an automorphism of VOSA $V(H,{\Bbb
Z}+\frac{1}{2})$ as follows: For any $a_1(-n_1)\cdots a_s(-n_s)\in
V(H,\Bbb Z+1/2),$ $$\tau(a_1(-n_1)a_2(-n_2)\cdots a_m(-n_m))=
(\tau a_1)(-n_1)(\tau a_2)(-n_2)\cdots (\tau a_m)(-n_m).$$ Let
$o(\tau \sigma)=N.$ We decompose $H$ into eigenspaces with respect to
the $\tau \sigma $ and $\tau $ as follows:
\begin{eqnarray}
& &H=\oplus_{r\in \Z/N\Z}H^{r*} \label{g10.9}\\
& & H=\oplus_{r\in \Z/N_0\Z}H^{r} \label{g10.10} 			
\end{eqnarray}
where $H^{r*}=\{v\in H|\tau\sigma v=e^{2\pi ir/N}v\},$ and $H^{r}=\{v\in H|\tau v=e^{2\pi ir/N_0}v\}.$

Let $l_0=\dim H^{0*}.$ As before we need to consider two separate cases:  
$l_0$ is even or odd.
If $l_0=2k_0$ for some nonnegative integer $k_0$, we have 
$$H^{0*}=\sum_{i=1}^{k_0}{\Bbb C} h_i+\sum_{i=1}^{k_0}{\Bbb C}h_i^*$$
with $(h_i,h_j)=(h_i^*,h_j^*)=0,  (h_i,h_j^*)=\delta_{i,j}.$

Let $l_r=\dim H^{r*}$ with $r\neq 0.$ If $r\neq N-r,$ we fix 
bases  
${b_{r,1},b_{r,2},\cdots b_{r,l_r}} \in H^{r*}$ and $b_{r,1}^{*},b_{r,2}^{*},\cdots b_{r,l_r}^{*} \in H^{(N-r)*}$
such that  $(b_{r,i},b_{r,j}^*)=(b_{r,j}^*,b_{r,i})=\delta_{i,j}.$
If $r=N-r,$ let  $\{c_1,c_2,\cdots c_{l_{\frac{N}{2}}}\}$
be an orthonormal basis of  $H^{\frac{N}{2}*}.$

Then $M=\bigotimes_{r=1}^{N-1}\Lambda[b(-n)|n\in \frac{r}{N},n>0,b\in
H^{r*}]\bigotimes \Lambda[h_i(-n),h_i^*(-m)|n,m \in {\Bbb Z},n>0,m \geq
0]$ is an irreducible $\tau$-twisted $V(H,\Z+\frac{1}{2})$-module so
that for $u\in H^{r*},$
\begin {eqnarray*}
Y_{M}(u(-\frac{1}{2}),z)=u(z)=\sum_{n\in \frac {r}{N}+{\Bbb Z}}u(n)z^{-n-1/2}
\end{eqnarray*}
(see \cite{Li2}).
Note that $b_{r,i}(n)$ acts as $\frac{\partial}{\partial b_{r,i}^{*}(-n)}$ if $n$ is positive and acts as multiplication by $b_{r,i}(n)$ if $n$ is negative.
Similarly for $b_{r,i}^*(n).$ The 
$c_i(n)$ acts as $\frac{\partial}{\partial c_i(-n)}$ if $n$ is positive and acts as multiplication by $c_i(n)$ if $n$ is negative.
 
 Also,  $h_i(n)$ acts as $\frac{\partial}{\partial h_i^*(-n)}$ if $n$ is 
 nonnegative, and acts as multiplication by $h_i(n)$ if $n$ is negative;
 $h_i^*(n)$ acts as $\frac{\partial}{\partial h_i(-n)}$ if $n$ is positive, 
and acts as multiplication by $h_i^*(n)$ if $n$ is nonnegative.

One can easily  calculate that  
$$\Omega(M)=\Lambda[h_i^{*}(0)|h_i^{*}\in H^{0*},i=1,2,\cdots k_0].$$
 So $\dim \Omega(M)=2^{k_0}.$

\begin{prop}\label {p7.3} If $\dim H^{0*}=l_0=2k_0$ then
 $M=\bigotimes_{r=1}^{N-1}\Lambda[b(-n)|n\in \frac{r}{N}+\Z,n>0,b\in H^{r*}]\bigotimes \Lambda[h_i(-n),h_i^*(-m)|n,m \in {\Bbb Z},n>0,m \geq 0]$
is the unique $\tau$-twisted irreducible $V(H,{\Bbb Z}+\frac{1}{2})$-module.
\end{prop}

\pf As in the proof of Proposition \ref{p7.1}, it is sufficient to
show that $\dim A_{\tau}(V(H,{\Bbb Z}+\frac{1}{2})) \leq
2^{l_0}=\dim \End\Omega(M).$

If $a\in H^{r*}$ with $r \neq 0,$ by Lemma  \ref{l3.2} (i) we see  that  for $m\geq0 , m \in {\Bbb Z}, b \in H$
\begin{eqnarray*}
 \Res_z \frac{(1+z)^{\frac{r}{N}-\frac{1}{2}}}{z^{1+m}}a(z)b=\sum_{l=0}^{\infty}\binom{\frac{r}{N}-\frac{1}{2}}{l}a(-m-\frac{1}{2}+l)b\in  O_{\tau}(V(H,{\Bbb Z}+\frac{1}{2})).
\end{eqnarray*}
So using the same calculation done in  Proposition \ref{p7.1}, 
we conclude that $A_{\tau}(V)$ is spanned by  $$h_1(-1/2)^{s_1}\cdots h_{k_0}(-1/2)^{s_{k_0}}h_1^*(-1/2)^{t_1}\cdots h_{k_0}^*(-1/2)^{t_{k_0}}$$ with $s_i,t_i=0,1.$ Hence $\dim A_{\tau}(V(H,{\Bbb Z}+\frac{1}{2})) \leq 2^{l_0},$
as desired. \qed  
\bigskip

If $\dim H^{0*}=l_0=2k_0+1,$ we can write $H^{0*}$ as follows:

$$H^{0*}=\sum_{i=1}^{k_0}{\Bbb C} h_i+\sum_{i=1}^{k_0}{\Bbb C}h_i^*+{\Bbb C}e$$
 with $(h_i,h_j)=(h_i^*,h_j^*)=0,  (h_i,h_j^*)=\delta_{i,j},(e,h_i)=(e,h_i^*)=0,(e,e)=2.$
 Let 
$$ W^0(H,{\Bbb Z})\!=\!\wedge[h_i(-n),h_i^{*}(-m),e(-n)|n\in\Z_{>0}, m\in \Z_{\geq0}]\!=\!W^0(H,{\Bbb Z})^{even}\oplus W^0(H,\Z)^{odd}.$$
Also define
$$V^0_{\pm}(H,\Z)=(1\pm e(0))W^0(H,{\Bbb Z})^{even}\oplus (1\mp e(0))W^0(H,{\Bbb Z})^{odd}.$$ 

Then  
 $$ M_{\pm}=V^0_{\pm}(H,\Z)\bigotimes \Lambda[a(-n)|a\in H^{r*},n\in   \frac{r}{N}+{\Bbb Z},1\leq r\leq N-1, n>0]$$
are irreducible $\tau$-twisted $V(H,{\Bbb Z}+\frac{1}{2})$- modules (see 
\cite{Li2}). 
The actions of $a(n),$ $h_i(m),$ $h_i^*(m)$ for 
$a \in H^{r*}$ for $r\neq 0$ are the same as before. The $e(n)$ acts as 
$2\frac{\partial}{\partial e(-n)}$ if $n>0$ and as multiplication
by $e(n)$ if $n\leq 0.$  The proof of Proposition \ref{p7.3} gives

\begin{prop}\label {p7.4} If $\dim H^{0*}=2k_0+1$ is odd, $V(H,{\Bbb Z}+\frac{1}{2})$ has exactly two $\tau$-twisted irreducible modules $M_{\pm}$ up to isomorphism.
\end{prop}

\end{document}